\documentclass[12pt]{article}
\usepackage[DIV12]{typearea} 
\usepackage[utf8]{inputenc}
\usepackage{amsmath}
\usepackage{amssymb}
\usepackage{amsfonts}
\usepackage{textcomp}
\usepackage{stmaryrd}
\usepackage{hyperref}
\usepackage[all]{xy}

\begin{document}
\newtheorem{thm}{Theorem}[section]
\newtheorem{prop}[thm]{Proposition}
\newtheorem{lem}[thm]{Lemma}
\newtheorem{defi}[thm]{Definition}
\newtheorem{defilemma}[thm]{Definition/Lemma}
\newtheorem{cor}[thm]{Corollary}
\newtheorem{exa}[thm]{Example}
\title{Moduli Spaces of Semistable Pairs in Donaldson$-$Thomas Theory}
\author{Malte Wandel\\\\Freie Universität Berlin, Leibniz Universität Hannover\\\\
current address:\\
Research Institute of Mathematical Sciences\\Kyoto University\\
Kitashirakawa Oiwake cho, Sakyo-ku\\
Kyoto 606-8502, JAPAN\\e-mail: wandel@kurims.kyoto-u.ac.jp\\
}

\maketitle
\begin{abstract}
Let $(X,\mathcal{O}_X(1))$ be a polarized smooth projective variety over the complex numbers. Fix $\mathcal{D}\in \mathrm{coh}(X)$ and a nonnegative rational polynomial $\delta$. Using GIT we contruct a coarse moduli space for $\delta$-semistable pairs $(\mathcal{E},\varphi)$ consisting of a coherent sheaf $\mathcal{E}$ and a homomorphism $\varphi\colon \mathcal{D}\rightarrow \mathcal{E}$. We prove a chamber structure result and establish a connection to the moduli space of coherent systems constructed by Le Potier in \cite{LeP} and \cite{LeP2}.\\\\
Keywords:\ moduli spaces, curve counting, stable pairs, coherent systems\\
MCS:\ 14D20, 14D22, 14N35
\end{abstract}
\tableofcontents
\setcounter{section}{-1}
\section{Introduction}
String theorists are highly interested in counting curves on Calabi$-$Yau threefolds. This can be done by integrating over a virtual cycle on the moduli space of these curves (cf.\ \cite{LT}). The arising moduli problems are not compact and there have been different approaches to their compactification (cf.\ \cite{PT}). Following Pandharipande and Thomas we consider pairs $(\mathcal{E}, s)$ consisting of a coherent sheaf $\mathcal{E}$ with one dimensional support and a section $s\in H^0(\mathcal{E})$. Such a pair is called stable if firstly the sheaf is pure and secondly $s$ considered as a homomorphism $\mathcal{O}_X\rightarrow\mathcal{E}$ is generically surjective. Thus a stable pair $(\mathcal{E}, s)$ provides us with a Cohen$-$Macaulay curve $C_\mathcal{E}=\mathrm{supp}\mathcal{E}$ and a finite number of points on this curve, namely the cokernel of $s$.\\
In \cite{LeP2} the moduli space for such stable pairs is contructed. More generally Le Potier considered so-called coherent systems $(\Gamma, \mathcal{E})$ consisting of a sheaf $\mathcal{E}$ together with a subspace $\Gamma\subseteq H^0(\mathcal{E})$.
Since the pairs introduced above can only cover the case of irreducible curves there is a need for generalizations of the notion of stable pairs. We want a section $s$ for every irreducible component of $C_\mathcal{E}$. Thus we should consider pairs $(\mathcal{E},\varphi)$ with a homomorphism $\varphi\colon\mathcal{O}_X^r\rightarrow \mathcal{E}$ or even more generally with a homomorphism $\varphi\colon \mathcal{D}\rightarrow \mathcal{E}$ for an arbitrary but fixed coherent sheaf $\mathcal{D}$. There is a generalized notion of stability for such pairs.\\
In this article we will construct a coarse moduli space for semistable pairs on an arbitrary polarized smooth projective variety and relate these moduli spaces to the moduli spaces of coherent systems in the case $\mathcal{D}=\mathcal{O}_X^r$.\\
We will now give a more detailed overview of the content of this article. In the first section we will define the generalized notion of stability depending on a parameter $\delta\in \mathbb{Q}[x]$ and discuss some basic properties of semistable pairs. Next in Section 2 we will prove the boundedness of the family of $\delta$-semistable pairs which will enable us to define the parameter space for our moduli problem in Section 3. The core of this article is contained in Section 4 where we perform the GIT construction of the moduli space. In Section 5 we prove the usual chamber structure result summarizing how the moduli spaces change if we vary the parameter $\delta$. Last but not least we show in Section 6 that the moduli space of coherent systems by Le Potier can be obtained from the moduli space of stable pairs as a quotient by a group action.\\
\\
Notations: The main guidelines for the notations used in this article are \cite{Har} and \cite{HL}. By a scheme one should always think of a scheme of finite type over some fixed algebraically closed field $k$ of characteristic zero. If $V$ is a finite dimensional $k$-vector space we let $\mathbb{P}(V):=(V\backslash\{0\})/\sim $ be the set of lines passing through the origin. For any vector $v\in V$ we denote its equivalence class in $\mathbb{P}(V)$ by $[v]$. Finally by PGL$_n$ we denote the quotient GL$_{n+1}/k^\star$Id. Thus PGL$_n$ is the automorphism group of $\mathbb{P}^n$.\\
\\
I want to thank Alexander Schmitt for refering to me such a challenging topic for a diploma thesis, for making it possible for me to travel to London and for supporting my work, especially the GIT part. Furthermore I have to thank Richard Thomas for fantastic hospitality, for giving me enlightening insights into modern geometry and of course for coming up with this moduli problem. At last I want to thank Klaus Altmann for co-supervising my thesis and R.W. for brushing up my English.\\
Finally I want to thank the referee for his suggestions to improve this note. This article is published in Manuscripta Mathematica (DOI: 10.1007/s00229-015-0729-7). The final publication is available at\\
\href{http://www.springer.com/-/1/45ec2ebe73104687a830eb127441ea9b}{http://www.springer.com/-/1/45ec2ebe73104687a830eb127441ea9b}.

\section{Semistable Pairs}
In order to define pairs we first fix some notation. Let $(X,\mathcal{O}_X(1))$ be a polarized smooth projective variety, $\mathcal{D}$ a coherent $\mathcal{O}_X$-module and $\delta$ a rational polynomial $\geq 0$, i.e., $\delta(m)\geq 0$ $\forall$ $m \gg 0.$

\begin{defi}
\em A pair $(\mathcal{E},\varphi)$ \em consists of a coherent $\mathcal{O}_X$-module $\mathcal{E}$ and a homomorphism $\varphi\colon\mathcal{D}\rightarrow \mathcal{E}$. By $P=P_ \mathcal{E}$ we denote the \em Hilbert polynomial \em of $\mathcal{E}$. This is a polynomial of degree $d=\dim\mathcal{E}$ and  its leading coefficient is just the rank of the sheaf $\mathcal{E}$ which we denote by $r$ or $r_\mathcal{E}$. We call a pair \em pure \em if $\mathcal{E}$ is pure. \em A homomorphism of pairs \em $\alpha\colon (\mathcal{E},\varphi)\rightarrow (\mathcal{E}',\psi)$ is a homomorphism of $\mathcal{O}_X$-modules $\alpha\colon \mathcal{E}\rightarrow \mathcal{E} '$ such that there is a scalar $\lambda\in k$ making the following diagramm commute
\[\xymatrix{ \mathcal{D}\ar[r]^{\lambda\cdot \mathrm{id}} \ar[d]_{\varphi} & \mathcal{D} \ar[d]^{\psi}\\
\mathcal{E} \ar[r]_{\alpha}& \mathcal{E}'&.}\]
In the obvious way we define the notion of \em an isomorphism of pairs\em.
\end{defi}

\begin{lem}\label{lemscmul}
Let $(\mathcal{E},\varphi)$ be a pair. Then for any $\lambda\in k^\star$ $(\mathcal{E},\lambda\varphi)$ is isomorphic to $(\mathcal{E},\varphi)$.
\end{lem}
\em Proof: \em Obvious. \hfill $\Box$

\begin{defi}
We define \em the Hilbert polynomial of a pair \em $(\mathcal{E},\varphi)$ to be
\[P_{(\mathcal{E},\varphi)} := P_{\mathcal{E}} + \epsilon (\varphi)\delta,\]
where $\epsilon(\varphi)=1$ if $\varphi\neq 0$ and $0$ otherwise. For any subsheaf $\mathcal{F}\subseteq \mathcal{E}$ we define \em the induced homomorphism \em $\varphi '$ to be equal to $\varphi$ if $\mathrm{im}\varphi\subseteq \mathcal{F}$ and $0$ otherwise. For the correspoding quotient $\mathcal{G}=\mathcal{E}/\mathcal{F}$ \em the induced homomorphism \em $\varphi ''$ is defined to be the composition of $\varphi$ with the quotient map. It is $0$ if and only if $\mathrm{im}\varphi\subseteq \mathcal{F}$. Thus we easily see that the Hilbert polynomial of pairs is additive in exact sequences.\\
\em The reduced Hilbert polynomial of a pair \em $(\mathcal{E},\varphi)$ is defined as
\[ p_{(\mathcal{E},\varphi)}= \frac{P_{(\mathcal{E},\varphi)}}{r_\mathcal{E}}.\]
\end{defi}

\begin{defi}
A pair $(\mathcal{E},\varphi)$ is called \em (semi)stable with respect to $\delta$ (or $\delta$-(semi)stable) \em if for every saturated submodule $\mathcal{F}\subseteq \mathcal{E}$ of rank $r'$ we have
\begin{eqnarray*}
P_{(\mathcal{F},\varphi ')} &(\leq)& r'p_{(\mathcal{E},\varphi)}.
\end{eqnarray*} 
\end{defi}
\em Remark\em : For $\delta=0$ this condition does not depend on $\varphi$ at all and we get the usual stability condition for sheaves. Thus the corresponding moduli-space is constructed and discussed in \cite{HL}. Therefore we will consider only strictly positive stability parameters $\delta$ from now on.\\
Similarly, if the homomorphism $\varphi$ is trivial the stability condition does not depend on $\delta$ and is again equivalent to the usual stability condition for sheaves. Thus during the construction of the moduli space of stable pairs we will restrict to pairs with nontrivial homomorphism $\varphi$.\\
\\
There is another good reason for not considering pairs with trivial homomorphism. Assume we have constructed a moduli space $M$ of $\delta$-semistable pairs allowing trivial homomorphisms. Fix an arbitrary semistable pair $(\mathcal{E},\varphi)$ with nontrivial $\varphi$ and consider the family $(\mathcal{E},\lambda\varphi)_{\lambda\in k}$ parametrized by the affine line. Assume that $(\mathcal{E},0)$ is as well semistable. This family would lead to a classifying map $f\colon\mathbb{A}^1\rightarrow M$. For every $\lambda\in k^\star$ we get the same point $f(\lambda)\in M$ (cf.\ \ref{lemscmul}). Thus for $\lambda=0$ we have to end up with the same point. But this time we have a point with trivial homomorphism. Thus we see that in such a case our moduli space would degenerate and parametrize sheaves only.

\begin{prop}\label{propssquot}
A pair is $\delta$-(semi)stable if and only if for every pure quotient $\mathcal{G}$ of $\mathcal{E}$ of rank $r''$ with induced homomorphism $\varphi ''$ we have:
\begin{eqnarray*}
P_{(\mathcal{G},\varphi '')}& (\geq)& r'' p_{(\mathcal{E},\varphi)}.
\end{eqnarray*}
\end{prop}
\em Proof: \em Let $\mathcal{G}$ be a pure quotient of $\mathcal{E}.$ Recall that by definition we have an exact sequence
\[0\rightarrow\mathcal{F}\rightarrow\mathcal{E}\rightarrow\mathcal{G}\rightarrow 0\]
with $\mathcal{F}$ a saturated submodule of $\mathcal{E}$. Also recall, that the rank and the Hilbert polynomial of pairs are additive. Now the claim follows just by substituting $P_\mathcal{G}$ and $r_\mathcal{G}$ by the appropriate expressions of $\mathcal{F}$. \hfill $\Box$

\begin{lem}\label{lemstabaut}
Let $\alpha\colon(\mathcal{E},\varphi)\rightarrow (\mathcal{E}',\psi)$ be a nontrivial homomorphism between two stable pairs of the same reduced Hilbert polynomial. Then $\alpha$ is an isomorphism.
\end{lem}
\em Proof: \em On $\mathcal{F}=\mathrm{im}\alpha$ there are two induced homomorphisms: the first coming from $\mathcal{E}$ called $\varphi ''$, the second induced by $\psi$ on $\mathcal{E}'$ called $\psi '$. It is easy to compute that $\psi '$ is trivial if and only if $\varphi ''$ is trivial. Thus the Hilbert polynomial of the pairs $(\mathcal{F},\varphi '')$ and $(\mathcal{F},\psi ')$ coincide. Since $(\mathcal{E},\varphi)$ and $(\mathcal{E}',\psi)$ have the same reduced Hilbert polynomial the stability assumption yields $\mathcal{F}=\mathcal{E}'$ and $\ker\alpha =0$. \hfill $\Box$

\begin{prop}\label{propjordan}
Let $(\mathcal{E},\varphi)$ be a $\delta$-semistable pair with reduced Hilbert polynomial $p$. There is a filtration
\[ 0 = \mathcal{E}_0 \subsetneq \mathcal{E}_1 \subsetneq \dots \subsetneq \mathcal{E}_r=\mathcal{E},\]
such that every factor $\mathrm{gr}_i(\mathcal{E})=\mathcal{E}_i/\mathcal{E}_{i-1}$ with its induced homomorphism is stable with reduced Hilbert polynomial $p$. Such a filtration is called \em a Jordan$-$Hölder filtration \em and the graded object $\mathrm{gr}(\mathcal{E}):=\oplus_i\mathrm{gr}_i(\mathcal{E})$ is independent of the choice of the filtration. It naturally inherits an induced homomorphism $\mathrm{gr}(\varphi)\colon \mathcal{D}\rightarrow \mathrm{gr}(\mathcal{E})$.
\end{prop}
\em Proof: \em If $(\mathcal{E},\varphi)$ is stable there is nothing to prove. Otherwise there is a subsheaf $\mathcal{F}\subseteq \mathcal{E}$ with induced homomorphism $\varphi '$ such that $p_{(\mathcal{F},\varphi ')}=p$. It is therefore semistable. Choose such a subsheaf $\mathcal{F}$ which is maximal with this property. Then the quotient $\mathcal{E}/\mathcal{F}$ is stable with reduced Hilbert polynomial $p$ by construction. We can now proceed in the same way with $\mathcal{F}$ constructing another subsheaf and so on. We end up with a filtration with stable factors which has to be finite since the rank is decreasing in every step. For a proof of the remaining assertions we refer to \cite{HL} Proposition 1.5.2 or \cite{HL2} Proposition 1.13. \hfill $\Box$\\
\\
\em Remark: \em One can easily see that for every semistable pair $(\mathcal{E},\varphi)$ with nontrivial $\varphi$ the homomorphism $\mathrm{gr}(\varphi)$ of the graded object is nontrivial and its image is contained in exactly one summand of $\mathrm{gr}(\mathcal{E})$.

\begin{defi}
Two $\delta$-semistable pairs $(\mathcal{E},\varphi)$ and $(\mathcal{E}',\psi)$ are called \em S-equivalent \em if their graded objects $(\mathrm{gr}(\mathcal{E}),\mathrm{gr}(\varphi))$ and $(\mathrm{gr}(\mathcal{E}'),\mathrm{gr}(\psi))$ are isomorphic.
\end{defi}
From now on until the end of this article we will only consider semistable pairs with nontrivial homomorphism $\varphi$. Thus whenever we say a pair is semistable, we mean semistable with nontrivial homomorphism. For later use we want to formulate a characterisation of the stability condition which is equivalent to the one given before if we restrict to the case of nontrivial homomorphisms. We first need another definition.
\begin{defi}
For every pair $(\mathcal{E},\varphi)$ and every exact sequence $0\rightarrow\mathcal{F}\rightarrow\mathcal{E}\rightarrow\mathcal{G}\rightarrow 0$ we define
\begin{eqnarray*}
 \epsilon (\mathcal{F})  := & \begin{cases} 1  & \text{if }\text{\em im\em}\varphi \subseteq \mathcal{F} \\ 0  & \text{otherwise}  \end{cases} & \text{ , and} \\
\epsilon (\mathcal{G}) := & 1- \epsilon (\mathcal{F}).
\end{eqnarray*}
\end{defi}

\begin{lem}
A pair $(\mathcal{E},\varphi)$ with $\varphi\neq 0$ is $\delta$-(semi)stable if and only if for every saturated subsheaf $\mathcal{F}\subseteq\mathcal{E}$ of rank $r'$ the following inequality holds:
\begin{eqnarray*}
P_\mathcal{F} + \epsilon(\mathcal{F})\delta& (\leq)& \frac{r'}{r}(P + \delta).
\end{eqnarray*}
\end{lem}
\em Proof: \em Obvious. \hfill $\Box$

\begin{prop}
Let $(\mathcal{E},\varphi)$ be a $\delta$-semistable pair. Then $\mathcal{E}$ is pure.
\end{prop}
\em Proof: \em Let $(\mathcal{E},\varphi)$ be a semistable pair and let $\mathcal{T}:=T_{d-1}(\mathcal{E})$ denote the maximal subsheaf of strictly smaller dimension. This is saturated and of rank zero. Thus the semistability condition yields
\begin{eqnarray*}
P_\mathcal{T} + \epsilon(\mathcal{T})\delta  \leq  \frac{r'}{r}(P + \delta)  =  0, 
\end{eqnarray*}
or equivalently $P_\mathcal{T}\leq - \epsilon(\mathcal{T})\delta \leq 0$, hence $\mathcal{T}=0$. \hfill $\Box$

\begin{prop}\label{propgensur}
If $\deg\delta\geq \dim X$ for every pair $(\mathcal{E},\varphi)$ the following two assumptions are equivalent:
\begin{itemize}
\item[\em(i)\em] $(\mathcal{E},\varphi)$ is $\delta$-semistable,
\item[\em(ii)\em] $\varphi$ is generically surjective and $\mathcal{E}$ is pure.
\end{itemize}
\end{prop}
\em Proof: \em (i) $\Rightarrow$ (ii): If $\varphi$ was not generically surjective, there would exist a saturated module $\mathcal{F}$ satisfying im$\varphi\subseteq \mathcal{F}\nsubseteq\mathcal{E}$ (for instance the saturation of im$\varphi$). Now semistability yields
\[ \frac{P_\mathcal{F}}{r_\mathcal{F}}+\frac{\delta}{r_\mathcal{F}}\leq \frac{P_\mathcal{E}}{r_\mathcal{E}}+\frac{\delta}{r_\mathcal{E}}. \]
For $\deg\delta\geq\deg P_\mathcal{E}$ it follows that $r_\mathcal{E} \leq r_\mathcal{F}$, thus we get $r_\mathcal{E} = r_\mathcal{F}$ because $\mathcal{F}\subseteq\mathcal{E}$ and we conclude $\mathcal{F}=\mathcal{E}$, because $\mathcal{F}$ was assumed to be saturated. This contradicts the assumption im$\varphi\subseteq \mathcal{F}\nsubseteq\mathcal{E}$. \\
(ii) $\Rightarrow$ (i): If $\varphi$ is generically surjective then for every saturated submodule $\mathcal{F}\subseteq\mathcal{E}$ we have im$\varphi\not\subseteq\mathcal{F}$, hence $\epsilon(\mathcal{F})=0$. Now if $(\mathcal{E},\varphi)$ was not semistable there would exist a saturated submodule $\mathcal{F}$ satisfying the destablizing condition:
\[P_\mathcal{F} > \frac{r'}{r}(P+\delta).\]
Since $\mathcal{E}$ is pure we have $r'\neq 0$. Thus for $\deg\delta>d$ we get $\delta<0$ and the same follows for $\deg\delta=d$ if we compare the leading coefficients of the polynomials. \hfill $\Box$\\
\\
\em Remark: \em By this proposition we see, that for large $\delta$ the moduli space we are interested in can be realized as some Quot-scheme parametrizing certain quotients of our fixed module $\mathcal{D}$. Thus from now on we will assume $\deg\delta<\dim X$.

\begin{prop}\label{propdeform}
If a pair $(\mathcal{E},\varphi)$ can be deformed into a pure pair, then there is a pure sheaf $\mathcal{H}$ and a morphism $\psi: \mathcal{E}\rightarrow \mathcal{H}$ satisfying $\ker\psi=\mathcal{T}(\mathcal{E})$. In particular if we set $\varphi_\mathcal{H}= \psi \circ \varphi$ we get a pure pair $(\mathcal{H},\varphi_\mathcal{H})$.
\end{prop}
\em Proof: \em The condition on $(\mathcal{E},\varphi)$ says there is a smooth connected curve $\mathcal{C}$ and a flat family $(\mathcal{E}_\mathcal{C},\varphi_\mathcal{C})$ on $\mathcal{C} \times X $ such that $(\mathcal{E}_0,\varphi_0) \cong (\mathcal{E},\varphi)$ for some closed point $0\in \mathcal{C}$ and $(\mathcal{E}_t,\varphi_t)$ pure for all $t\neq 0$. In particular $\mathcal{E}$ deforms into a pure sheaf. Now our claim follows from \cite{HL}, Proposition 4.4.2. \hfill $\Box$\\
\\
\em Remark \em : Note that a priori $\varphi_\mathcal{H}$ could be trivial. This might occur if im$\varphi$ is contained in $T_{d-1}(\mathcal{E})\subseteq\mathcal{E}$. Therefore in order to show that such a $(\mathcal{H},\varphi_\mathcal{H})$ is semistable we should first show that $\varphi_\mathcal{H}$ is infact nontrivial.

\section{Boundedness}
In this section we prove that the set of $\delta$-semistable pairs with fixed Hilbert polynomial is bounded. Furthermore we deduce an important stability criterion which we will need for the construction of the moduli space in Section 4.

\begin{prop}\label{propbound}
Let $P$ and $\delta\geq 0$ be polynomials. Then there is a constant $C$ depending only on $P$ and $\mathcal{D}$ such that for every $\mathcal{O}_X$-module $\mathcal{E}$ occuring in a $\delta$-semistable pair we have $\mu_{\mathrm{max}}(\mathcal{E}) \leq C$. In particular, the family of pairs which are semistable with respect to any stability parameter $\delta$ having the fixed Hilbert polynomial $P$ is bounded.
\end{prop}
\em Remark: \em Note that the uniform bound is independent of $\delta$.\\
\em Proof: \em Let $\mu_P$ denote the slope of $P$ and $\delta_1$ denote the coefficient of $\delta$ in degree $d-1$. It is $\delta_1\geq 0$. If $\mathcal{F}$ is a submodule of $\mathcal{E}$ of rank $r'$ satisfying im$\varphi \subseteq \mathcal{F}$ then the semistability condition yields
\[
\mu_\mathcal{F} + \frac{\delta_1}{r'}\leq \mu_P + \frac{\delta_1}{r}.
\]
And since $r'\leq r$ we have $\mu_\mathcal{F}\leq \mu_P$. Now let $\mathcal{F}\subseteq\mathcal{E}$ be an arbitrary submodule. We have an exact sequence \[0\rightarrow\mathcal{F}\rightarrow\mathcal{F} + \text{im}\varphi \rightarrow \mathcal{G}\rightarrow 0,\] where $\mathcal{G}=\text{im}\varphi/(\mathcal{F}\cap\text{im}\varphi)$ is a quotient of im$\varphi$, so a fortiori a quotient of $\mathcal{D}$. Set $\mathcal{H} := \mathcal{F} + \mathrm{im}\varphi$ and note that $\mathcal{H}$ contains im$\varphi$, so we can apply the first part of the proof. By the additivity of the degree we have:
\begin{eqnarray*}
\mu_\mathcal{F} &=& \frac{\mu_{\mathcal{H}}r_{\mathcal{H}}-\mu_\mathcal{G}r_\mathcal{G}}{r_\mathcal{F}} \\
&\leq & \frac{\mu_Pr_\mathcal{H}-\mu_{\text{min}}(\mathcal{D})r_\mathcal{G}}{r_\mathcal{F}}\\
&\leq & \mu_Pr-\mu_{\text{min}}(\mathcal{D})\frac{r_\mathcal{G}}{r_\mathcal{F}}.
\end{eqnarray*}
Now depending on the sign of $\mu_{\text{min}}(\mathcal{D})$ we find a uniform bound C for $\mu_\mathcal{F}$ by setting $C:=\max\{\mu_P,\mu_Pr-\mu_{\text{min}}(\mathcal{D})r,\mu_Pr-\mu_{\text{min}}(\mathcal{D})\frac{1}{r}\}$ \hfill. $\Box$\\
\\
Next we need a result to get an estimate for the global sections of certain sheaves. It follows from a more general result due to Simpson (cf.\ \cite{Sim}, Lemma 1.5). Therefore we define $\alpha_i$ to be the i-th coefficient of $P_{\mathcal{O}_X}$, the Hilbert polynomial of the structure sheaf. For any $\mathcal{O}_X$-module $\mathcal{E}$ we define 
\[
\widehat\mu(\mathcal{E}):= \frac{\mu_\mathcal{E}+\alpha_{d-1}}{\alpha_{d}}.\]
Note that $\widehat\mu(\mathcal{E})$ is just the coefficient of $P_\mathcal{E}/r$ in degree $d-1$.

\begin{prop}\label{propsimpsest}
Let $\mathcal{E}$ be a coherent $\mathcal{O}_X$-module of rank $r$ and dimension $d$ and let $C:=r(r+d)/2$. Then
\[
\frac{h^0(\mathcal{E}(m))}{r}\leq \frac{r-1}{r}\cdot \frac{1}{d!}[\widehat\mu_{\mathrm{max}}(\mathcal{E}) + C - 1 + m]^d_+
 + \frac{1}{r}\cdot \frac{1}{d!}[\widehat\mu(\mathcal{E}) +C -1 +m]^d_+.
\]
\end{prop}
\em Proof: \em \cite{HL}, Corollary 3.3.8. \hfill $\Box$

\begin{cor}
If in addition $\mathcal{E}$ is semistable we have
\[
\frac{h^0(\mathcal{E}(m))}{r}\leq \frac{1}{d!}[\widehat\mu(\mathcal{E}) +C -1 +m]^d_+.
\]
\end{cor}
\em Proof: \em Since $\mathcal{E}$ is semistable we have $\mu_{max}(\mathcal{E})\leq \mu(\mathcal{E})$.\hfill $\Box$\\
\\
Now we're able to proof a very important stability criterion which we will need for the GIT construction of the moduli space. It corresponds logically to \cite{HL}, Theorem 4.4.1 and we will follow their proof closely. First we fix a stability parameter $\delta$. Let $m_1$ be an integer such that $\delta(m)\geq 0$ $\forall m\geq m_1$.

\begin{prop}\label{propssrephrase}
There is an integer $m_0$ such that for every $m\geq m_0$ and every pure pair $(\mathcal{E},\varphi)$ with Hilbert polynomial $P$ and rank $r$ the following three assertions are equivalent:
\begin{itemize}
\item[\em (i)\em] $(\mathcal{E},\varphi)$ is $\delta$-(semi)stable,
\item[\em (ii)\em] $P(m)  \leq h^0(\mathcal{E}(m))$ and for all subsheaves $\mathcal{F}\subseteq\mathcal{E}$ of rank $0<r'<r$ we have
\begin{eqnarray*}
h^0(\mathcal{F}(m)) + \epsilon(\mathcal{F})\delta(m) & (\leq) &\frac{r'}{r}\big(P(m)+\delta(m)\big),
\end{eqnarray*} 
\item [\em (iii)\em] for all quotients $\mathcal{E}\rightarrow\mathcal{G}$ of rank $0<r''<r$ we have
\begin{eqnarray*}
\frac{r''}{r}\big(P(m)+\delta(m)\big)&(\leq)& h^0(\mathcal{G}(m)) +  \epsilon(\mathcal{G})\delta(m).
\end{eqnarray*} 
\end{itemize}
\end{prop}
\em Remark: \em Note that we may assume $m \geq m_1$.\\
\em Proof: \em (i) $\Rightarrow$ (ii): Since the family of semistable pairs $(\mathcal{E},\varphi)$ is bounded, there is an integer $m$ such that for any $\mathcal{E}$ occuring in such a pair we have $P(m)=h^0(\mathcal{E}(m))$. Now let $\mathcal{F}\subseteq\mathcal{E}$ be an arbitrary subsheaf of rank $r'$. To show (2) we may assume $\mathcal{F}$ to be saturated. First we assume that $\max\{\mu_P,\mu_Pr-\mu_{\text{min}}(\mathcal{D})\}=\mu_P=\mu_\mathcal{E}$ and distinguish two cases:
\begin{itemize}
\item[(1)] $\widehat\mu(\mathcal{F}) \geq \widehat\mu(\mathcal{E}) - C\cdot r-\delta_1$,
\item[(2)]$\widehat\mu(\mathcal{F}) < \widehat\mu(\mathcal{E}) - C\cdot r-\delta_1$,
\end{itemize}
where $C:=r(r+d)/2$ as before. In the family of subsheaves of type (1) $\widehat\mu(\mathcal{F})$ is bounded from below. Since $\widehat\mu(\mathcal{F})= (\mu_\mathcal{F}+\alpha_{d-1})/\alpha_{d}$ we easily find that $\mu_\mathcal{F}$ is aswell bounded from below. Since we are talking about saturated subsheaves only by Grothendieck's proposition (cf.\ Proposition \ref{propssrephrase}) the family of subsheaves of type (1) is bounded. Thus the set of Hilbert polynomials of this family is finite and by taking $m$ sufficiently large we have $h^0(\mathcal{F}(m))=P(\mathcal{F}(m))$ and
\begin{eqnarray*}
P(\mathcal{F}(m)) + \epsilon(\mathcal{F})\delta(m) (\leq) \frac{r'}{r}[P(m)+\delta(m)] 
& \Leftrightarrow & P(\mathcal{F}) + \epsilon(\mathcal{F})\delta (\leq) \frac{r'}{r}[P+\delta].
\end{eqnarray*}
Now we consider subsheaves of type (2). Note that $\widehat\mu_{\text{max}}(\mathcal{F})\leq \widehat\mu_{\mathrm{max}}(\mathcal{E})\leq \mu(\mathcal{E})$ by the assumption made at the beginning. Set $C':=r'(r'+d)/2$ and note that it is $>0$. By the estimate of Proposition \ref{propsimpsest} we have
\begin{eqnarray*}
\frac{h^0(\mathcal{F}(m))}{r'}\leq\\
 \frac{r'-1}{r'}\cdot\frac{1}{d!}[\widehat\mu_{\mathrm{max}}(\mathcal{F}) + C' - 
1 + m]^d_+ + \frac{1}{r'}\cdot \frac{1}{d!}[\widehat\mu(\mathcal{F}) +C' -1 +m]^d_+
\leq \\
\frac{r'-1}{r'}\cdot\frac{1}{d!}[\widehat\mu(\mathcal{E}) + C' - 
1 + m]^d_+ + \frac{1}{r'}\cdot \frac{1}{d!}[\widehat\mu(\mathcal{E}) +C'-C\cdot r
-\delta_1 -1 +m]^d_+ =\\
  \frac{m^d}{d!}+ \frac{m^{d-1}}{(d-1)!} \Big(\,
\underbrace{\widehat\mu(\mathcal{E}) + C' - 1 + \frac{1}{r'}(C\cdot r-\delta_1)}_{=:A} \Big)\, + O(d-2)
\end{eqnarray*}
where $O(d-2)$ stands for polynomials in $m$ of degree $\leq d-2$. Since $r'\leq r$ we have
\begin{eqnarray*}
A=\widehat\mu(\mathcal{E})-1+C'- \frac{r}{r'}C-\frac{\delta_1}{r'}<\widehat\mu(\mathcal{E})-1-\frac{\delta_1}{r'}.
\end{eqnarray*}
Now as we have said before $\widehat\mu(\mathcal{E})$ is the coefficient of $P/r$ in degree $d-1$. Thus for sufficiently large $m$ we get
\[
\frac{1}{r'}\big(h^0(\mathcal{F}(m))+\epsilon(\mathcal{F})\delta(m)\big) < \frac{m^d}{d!}+\frac{m^{d-1}}{(d-1)!} (\widehat\mu(\mathcal{E}) -1) < \frac{P(m)}{r} < \frac{P(m) + \delta (m)}{r}
\]
Now if $\max\{\mu_P,\mu_Pr-\mu_{\text{min}}(\mathcal{D})\}=\mu_Pr-\mu_{\text{min}}(\mathcal{D})$ we can run through all the arguments before substituting $\widehat\mu(\mathcal{E})$ by
\[\widehat\mu(\mathcal{E})r-\frac{\mu_{min}(\mathcal{D})+(r-1)\alpha_{d-1}}{\alpha_d}\] in the inequalities of the distinction of the cases.\\
(ii) $\Rightarrow$ (iii): For any quotient $\mathcal{E}\rightarrow\mathcal{G}$ we let $\mathcal{F}\subseteq\mathcal{E}$ denote the corresponding kernel. By (ii) we get:
\begin{eqnarray*}
h^0(\mathcal{G}(m)) + \epsilon(\mathcal{G})\delta(m)& \geq &h^0(\mathcal{E}(m)) -h^0(\mathcal{F}(m)) + \delta(m) - \epsilon(\mathcal{F})\delta(m) \\
&(\geq) &\frac{1}{r}\Big(\,rP(m)-r'P(m)+r\delta(m)-r'\delta(m)\Big)\,\\
& = & \frac{r''}{r}\Big(\,P(m)+\delta(m)\Big)\,.
\end{eqnarray*}
(iii) $\Rightarrow$ (i): Let $\mathcal{E}$ be an arbitrary sheaf satisfying (iii). We denote by $\mathcal{G}_{\mathrm{min}}$ the minimal destabilizing quotient sheaf of $\mathcal{E}$. Since it is semistable it has the nice property that $\mu_{\mathrm{min}}(\mathcal{G}_{\mathrm{min}})=\mu(\mathcal{G}_{\mathrm{min}})=\mu(\mathcal{G}_{max})$. By (iii) and Proposition \ref{propsimpsest} we have
\[
\frac{P(m)+\delta(m)}{r}-\frac{\epsilon(\mathcal{G}_{\mathrm{min}})\delta(m)}{r''} \leq\frac{h^0(\mathcal{G}_{\mathrm{min}}(m))}{r''} \leq \frac{1}{d!}[\widehat\mu(\mathcal{G}_{\mathrm{min}})+C-1+m]^d_+.
\]
Thus $\widehat\mu_{\mathrm{min}}(\mathcal{E})=\widehat\mu(\mathcal{G}_{\mathrm{min}})$ is bounded from below which is equivalent to $\widehat\mu_{\mathrm{max}}(\mathcal{E})$ being bounded from above. Hence we find the family of sheaves $\mathcal{E}$ satisfying (iii) to be bounded. Now we want to apply Proposition \ref{propssquot}. Thus let $\mathcal{G}$ denote an arbitrary quotient. We have either $\widehat\mu(\mathcal{G}) > \widehat\mu(\mathcal{E})+\delta_1/r$ but then we get a strict inequality in the semistability condition or $\widehat\mu(\mathcal{G}) \leq \widehat\mu(\mathcal{E})+\delta_1/r$, i.e., $\widehat\mu(\mathcal{G})$ is bounded from above and we can once more apply Grothendieck's result to get the boundedness of the family of the quotients in question. Again for large $m$ we have $h^0(\mathcal{G}(m))=P(\mathcal{G}(m))$ and 
\begin{eqnarray*}
P(\mathcal{G}(m)) + \epsilon(\mathcal{G})\delta(m) (\geq) \frac{r''}{r}[P(m)+\delta(m)] 
& \Leftrightarrow & P(\mathcal{G}) + \epsilon(\mathcal{G})\delta (\geq) \frac{r''}{r}[P+\delta].
\end{eqnarray*}
This finally shows (iii) $\Rightarrow$ (i) and finishes the proof.\hfill $\Box$\\
\em Remark: \em The proof shows that in (ii) and (iii) equality holds if and only if the subsheaf or the quotient, resp., is destabilizing.

\section{The Moduli Problem}
We will now define the moduli functor for stable pairs and present the parameter space which is somewhat more complicated than just a quot-scheme. Furthermore we discuss the natural linearizations on the parameter space and state the central theorems of this article, the existence of the moduli space of stable pairs.\\
Throughout this section we will fix a stability parameter $\delta$, a rational polynomial $P$ and a polarized smooth projective scheme $(X,\mathcal{O}_X(1))$.
\begin{defi}
We define a functor
\[
\mathcal{M}_{X,\delta}(\mathcal{D},P)\colon (\mathrm{Sch}/k)\text{\textdegree}\rightarrow (\mathrm{Sets})
\]
as follows: Let $S$ be a $k$-scheme of finite type. Define $\mathcal{M}_{X,\delta}(\mathcal{D},P)(S)$ to be the set of isomorphism classes of pairs $(\mathcal{E},\varphi)$ consisting of a coherent $S$-flat sheaf $\mathcal{E}$ on $X\times S$ (coming with the projections $\pi_X$ onto $X$ and $\pi_S$ onto $S$) and a homomorphism $\pi_X^\star\mathcal{D}\rightarrow \mathcal{E}$, such that for every closed point $s\in S$ the pair $(\mathcal{E}_s,\varphi|_{\pi_X^\star\mathcal{D}_s})$ is $\delta$-semistable with Hilbert polynomial $P$. For every morphism of $k$-schemes $f\colon S'\rightarrow S$ we obtain a map $\mathcal{M}_{X,\delta}(\mathcal{D},P)(f)\colon \mathcal{M}_{X,\delta}(\mathcal{D},P)(S)\rightarrow\mathcal{M}_{X,\delta}(\mathcal{D},P)(S')$ by pulling back $\mathcal{E}$ and $\varphi$ by $\mathrm{Id}_X\times f$. We define a subfunctor $\mathcal{M}^s_{X,\delta}(\mathcal{D},P)$ concerning only stable pairs.
\end{defi}

\begin{defi}
A scheme $M_{X,\delta}(\mathcal{D},P)$ is called \em a coarse moduli space of $\delta$-semistable pairs\em , if it correpresents the functor $\mathcal{M}_{X,\delta}(\mathcal{D},P)$.
\end{defi}
As we have seen before the family of isomorphism classes of sheaves with Hilbert polynomial $P$ occuring in a $\delta$-semistable pair is bounded. Therefore there is an $m_2$ such that for all integers $m\geq m_2$ the sheaves $\mathcal{E}$ are $m$-regular. Now we fix one $m\geq m_2,m_1,m_0$ (for notations see Proposition \ref{propssrephrase}) until the end of Section 4 and note that we may assume that $\mathcal{D}$ is $m$-regular as well.\\
Set $V=k^{P(m)}$. Let $Q :=\mathrm{Quot}_X(V\otimes\mathcal{O}_X(-m),P)$ be Grothendieck's Quot-scheme parametrizing quotients $q: V\otimes\mathcal{O}_X(-m) \rightarrow \mathcal{E}$ with Hilbert polynomial $P$. Recall that for large $l$ we have the very ample line bundles $\mathcal{L}_l$ on $Q$ which come from pullbacks of the universal bundle on $Q\times X$. Furthermore let $N:=\mathbb{P}\big(\mathrm{Hom}(H^0(\mathcal{D}(m)),V)\big)$ be the space of morphisms $ a\colon H^0(\mathcal{D}(m))\rightarrow V$ which is polarized by $\mathcal{O}_N(1)$.

\begin{lem}\label{lemassaq}
Let $(\mathcal{E},\varphi)$ be a $\delta$-semistable pair with Hilbert polynomial $P$. To this pair we can then associate in a natural way a pair $(a,q)\in N\times Q$ such that the induced map $H^0(q(m))$ is an isomorphism and $q\circ a = \varphi\circ \mathrm{ev}$:
\[ \xymatrix{H^0(\mathcal{D}(m))\otimes \mathcal{O}_X(-m) \ar[d]_a \ar[r]^/2em/{\mathrm{ev}} &
\mathcal{D} \ar[d]^\varphi \\
V\otimes\mathcal{O}_X(-m) \ar[r]_/1em/q  & \mathcal{E}&.
} \]
Here `\em ev\em' denotes the natural evaluation map.
\end{lem}
\em Proof: \em Since $\mathcal{E}$ is $m$-regular $\mathcal{E}(m)$ is globally generated. Since the Hilbert polynomial of $\mathcal{E}$ is P, we have $h^0(\mathcal{E}(m))=P(m)=\dim V$. If we choose any isomorphism $V\rightarrow H^0(\mathcal{E}(m))$ we find a surjection
\[ \rho\colon V\otimes \mathcal{O}_X\rightarrow \mathcal{E}(m).\]
Since taking tensor products is a right exact functor we obtain by tensoring with $\mathcal{O}_X(-m)$ another surjection
\[\rho(m)\colon V\otimes \mathcal{O}_X(-m)\rightarrow \mathcal{E}\]
and hence a point $q$ in $Q$. On the other hand we consider $\varphi\colon \mathcal{D}\rightarrow \mathcal{E}$. Tensoring with $\mathcal{O}_X(m)$, applying the global section functor and again choosing some isomorphism $V\rightarrow H^0(\mathcal{E}(m))$ we obtain a map
\[a:=H^0(\varphi(m))\colon H^0(\mathcal{D}(m))\rightarrow V.\]
Note that since $(\mathcal{E},\varphi)$ is assumed to be semistable the homomorphism $\varphi$ is not zero. Thus $a$ is not zero as well. Furthermore Lemma \ref{lemscmul} says that $(\mathcal{E},\varphi)$ and $(\mathcal{E},\lambda\varphi)$ are isomorphic pairs, thus $a$ is a well-defined point in $N$. It is clear from the construction that $q\circ a = \varphi\circ \mathrm{ev}$ and $H^0(q(m))$ is an isomorphism.\hfill $\Box$

\begin{prop}\label{propclsub}
There is a closed subscheme $\mathcal{Z}'\subseteq N\times Q$ such that for every pair
$(a,q)\in N\times Q$ the composition $q\circ a$ factors through the evaluation map \em ev \em (and thus induces a homomorphism $\varphi$) if and only if $(a,q)\in  \mathcal{Z}'$:
\[ \xymatrix{H^0(\mathcal{D}(m))\otimes \mathcal{O}_X(-m) \ar[d]^a \ar[r]^/2em/{\mathrm{ ev}} &
\mathcal{D} \ar@/^/@{ .>}[ldd]^\varphi \\
V\otimes\mathcal{O}_X(-m) \ar[d]^q  \\
\mathcal{E}&&.
} \]
\end{prop}
\em Proof: \em  We'll have a look at the relative version of the upper diagramm on $N\times Q\times X$:
\[\xymatrix{ 0 \ar[r] & \mathcal{K} \ar[r] \ar[rd]_f & H^0(\mathcal{D}(m))\otimes \pi ^\star_X \mathcal{O}_X(-m) \ar[r]^/2em/{\mathrm{ev}} \ar[d]^{\tilde{q} \circ \tilde a} & \pi ^\star_X\mathcal{D} \ar@{ .>}[ld]^{\tilde \varphi} \ar[r] & 0 \\
& & \pi_Q^\star\mathcal{E}&&.} \]
Here $\tilde q$ denotes the pullback of the universal quotient on $Q\times X$ and $\tilde a$ denotes the pullback of the universal homomorphism on $N$.\\
Now  $\tilde q \circ \tilde a$ factors through ev if and only if $f$ vanishes. Indeed, the if direction follows from the fundamental theorem on homomorphisms, the only if direction is an immediate consequence of the exactness of the upper row of the diagramm. Now the claim follows from \cite{Schm}, Proposition 2.3.5.1 and \cite{GS}, Lemma 3.1.\hfill $\Box$ \\
\\
The set of points $(a,q)$ where $\mathcal{E} = q(V\otimes \mathcal{O}_X(-m))$ is pure forms an open subset $\mathcal{U}\subseteq \mathcal{Z}'$. We let $\mathcal{Z}$ be its closure in $\mathcal{Z}'$.
Now SL$(V)$ acts diagonally on $\mathcal{Z}$:
\[ g\cdot (a,q)= (g \circ a, q \circ g^{-1}(-m)).
\]
Now $\mathcal{L}_l$ and  $\mathcal{O}_N(1)$ inherit natural linearizations of this SL$(V)$-action. Thus for integers $ n_1, n_2 $ the very ample line bundles $\mathcal{O}_\mathcal{Z}(n_1,n_2):=\pi ^\star_Q\mathcal{L}_l^{\otimes n_1}\otimes \pi ^\star_N \mathcal{O}_N(n_2)$ as well inherit natural SL$(V)$-linearizations. For some fixed $l$ we will choose $n_1,n_2$ to satisfy
\[\frac{n_2}{n_1} = \frac{P(l)\cdot\delta(m)-\delta(l)\cdot P(m)}{P(m)+\delta(m)}.\]

\begin{defilemma}\label{defilemr}
We define $\mathcal{R}\subseteq \mathcal{Z}$ to be the subset consisting of points $(a,q)$ corresponding to pairs which are $\delta$-semistable and the induced map $H^0(q(m))$ is an isomorphism. This is an open and \em SL\em$(V)$-invariant subset of $\mathcal{Z}$. Furthermore there is an open subset $\mathcal{R}^s\subseteq \mathcal{R}$ corresponding to $\delta$-stable pairs.
\end{defilemma}
\em Proof: \em Consider the universal homomorphism $\tilde{q}\circ \tilde{a}\colon H^0(\mathcal{D}(m))\otimes \pi_X^\star\mathcal{O}_X(-m)\rightarrow \pi_Q^\star\mathcal{E}$ on $\mathcal{Z}\times X$ and the projection $\pi\colon \mathcal{Z} \times X\rightarrow\mathcal{Z}$ which is a projective morphism since $X$ is projective. Let us consider the set $A$ of polynomials $P'$ occuring as a Hilbert polynomial of a quotient which destabilizes a pair corresponding to a point $(a,q)\in \mathcal{Z}$. For every such polynomial $P'$ its slope $\mu '$ satisfies
\[ \mu ' \leq \mu_P + \frac{\delta}{r_P}.\]
By \cite{Gro}, 2.5 the set $A$ is finite. For every $P'$ we consider the relative quot-scheme  $Q(P'):=\mathrm{Quot}_{\mathcal{Z}\times X/\mathcal{Z}}(\pi_Q^\star\mathcal{E},P')$. It comes with a morphism $f_{P'}\colon Q(P')\rightarrow \mathcal{Z}$. We denote the image of $f$ by $\mathcal{Z}(P')$. It is a closed subset of $\mathcal{Z}$. Now it is clear that a point $(a,q)\in\mathcal{Z}$ corresponds to a semistabe pair if and only if it is not in the finite and therefore closed union $\bigcup_{P'\in A} \mathcal{Z}(P')$.\hfill $\Box$

\begin{prop}
If a scheme $M$ is a categorical quotient for the \em SL\em$(V)$-action on $\mathcal{R}$, then it correpresents the functor $\mathcal{M}_{X,\delta}(\mathcal{D},P)$.
\end{prop}
\em Proof: \em Let $M$ be such a categorical quotient of $\mathcal{R}$. We have to show that for any $k$-scheme of finite type $S$ there is a map $\mathcal{M}_{X,\delta}(\mathcal{D},P)(S)\rightarrow \mathrm{Hom}(S,M)$. Therefore let $S$ be such a $k$-scheme and $(\mathcal{E},\varphi)$ be an element of $\mathcal{M}_{X,\delta}(\mathcal{D},P)(S)$. Thus 
\[V_\mathcal{E}:=p_\star(\mathcal{E}\otimes q^\star\mathcal{O}_X(m))\]
is a locally free sheaf on $S$ of rank $P(m)=\dim V$. We obtain a quotient on $X\times S$:
\[ \varrho_\mathcal{E}\colon p^\star V_\mathcal{E}\otimes q^\star\mathcal{O}_X(-m) \rightarrow \mathcal{E}.\]
Let $F_\mathcal{E}:=\mathbb{I}$som$(V\otimes \mathcal{O}_S,V_\mathcal{E})$ be the frame bundle of $V_\mathcal{E}$ together with the natural projection $\pi:F_\mathcal{E}\rightarrow S$ and the universal trivialisation of $V_\mathcal{E}$
\[f\colon V\otimes\mathcal{O}_{F_\mathcal{E}}\rightarrow \pi^\star V_\mathcal{E}.\]
Now from the quotient
\[ q_\mathcal{E}:=(id_X\times\pi)^\star\varrho_\mathcal{E} \circ\pi_{F_\mathcal{E}}^\star f\colon \mathcal{O}_X(-m)\otimes V\otimes\mathcal{O}_{F_\mathcal{E}}\rightarrow \pi_X\mathcal{E}\]
on $X\times F_\mathcal{E}$ we obtain its classifying morphism
\[\Phi_\mathcal{E}\colon F_\mathcal{E} \rightarrow Q. \]
On the other hand Lemma \ref{lemassaq} states that the homomorphism $\varphi$ yields a homomorphism of vector spaces $a\colon H^0(\mathcal{D}(m))\rightarrow V$. Lifting this map to GL$(V)$ and composing with the universal automorphism $\tau\colon \mathcal{O}_{\mathrm{GL}(V)}\otimes V \rightarrow\mathcal{O}_{\mathrm{GL}(V)}\otimes V$ we obtain a homomorphism on GL$(V)$
\[\mathcal{O}_{\mathrm{GL}(V)}\otimes H^0(\mathcal{D}(m))\rightarrow \mathcal{O}_{\mathrm{GL}(V)}\otimes V.\]
Again this yields a classifying map
\[\Xi_a\colon \mathrm{GL}(V)\rightarrow N.\]
Clearly there is a natural action of GL$(V)$ on $F_\mathcal{E}$ and GL$(V)$ itself. Moreover it is not hard to see that $S$ is a categorical quotient for the GL$(V)$-action on the product GL$(V)\times F_\mathcal{E}$. Alltogether we get a diagramm where vertical arrows are the maps from the categorical quotients:
\[\xymatrix@C=50pt{
\mathrm{GL}(V)\times F_\mathcal{E} \ar[r]^/0.5em/{\Xi_a\times\Phi_\mathcal{E}} \ar[d] & N\times Q \ar[d]\\
S \ar@{ .>}[r] & M
}\]
It follows from the construction that the map $\Xi_a\times\Phi_\mathcal{E}$ is GL$(V)$-equivariant and therefore yields a map between the categorical quotients $S\rightarrow M$.\hfill $\Box$\\
\\
The following result forms the center of this article. The proof will take all of Section 4.

\begin{thm}\label{thmconstr}
For sufficiently large $l$ the subset of points in the closure $\bar{\mathcal{R}}$ of $\mathcal{R}$ which are (semi)stable with respect to the $\mathrm{SL}(V)$-linearization coincides with the subset of points corresponding to $\delta$-(semi)stable pairs.
\end{thm}

\begin{thm}\label{thmmodspex}
Let $(X,\mathcal{O}_X(1))$ be a polarized smooth projective variety, $\mathcal{D}$ a coherent $\mathcal{O}_X$-module and $\delta$ a nonegative rational polynomial. Then there exists a coarse moduli space $M_{X,\delta}(\mathcal{D},P)$ of $\delta$-semistable pairs. Two pairs correspond to the same point in $M_{X,\delta}(\mathcal{D},P)$ if and only if they are S-equivalent. Moreover there is an open subset $M^s_{X,\delta}(\mathcal{D},P)\subseteq M_{X,\delta}(\mathcal{D},P)$ corresponding to stable pairs. It is a fine moduli space of $\delta$-stable pairs, i.e., it represents the functor $\mathcal{M}^s_{X,\delta}(\mathcal{D},P)$.
\end{thm}
\em Proof: \em The existence of the moduli space follows easily from Theorem \ref{thmconstr} and \cite{GIT}, Theorem 1.10 or \cite{Schm}, Theorem 1.4.3.8. The first remaining thing is to show that the coarse moduli space in fact parametrizes S-equivalence classes. Similarly as in the proofs of Proposition 3.3 in \cite{HL2} and Lemma 4.1.2 in \cite{HL} one can show that the orbit of any semistable pair $(\mathcal{E},\varphi)$ corresponding to a point in $\mathcal{R}^{\mathrm{ss}}$ also contains the graded object $($gr$(\mathcal{E}),$gr$(\varphi))$ and the orbits of these graded objects are closed. Thus S-equivalent pairs are mapped to the same point and since a geometric quotient map separates closed orbits we are done.\\
The second and last open statement is the fact that $M^s:=M^s_{X,\delta}(\mathcal{D},P)$ is indeed a fine moduli space. This is equivalent to the existence of a universal family on $X\times M^s$. Following Section 4.6 in \cite{HL} (in particular Proposition 4.6.2) the only thing we have to do is to show that there is a line bundle on $\mathcal{R}^s$ on which the center $Z:=k^\star\cdot\mathrm{Id}\subseteq\mathrm{GL}(V)$ acts with weight $1$. Such a line bundle can be named explicitly. Just remember that $\mathcal{R}$ was a subset of the product $N\times Q=\mathrm{Quot}_X(V\otimes\mathcal{O}_X(-m),P)\times
\mathbb{P}\big(\mathrm{Hom}(H^0(\mathcal{D}(m)),V)\big)$. Now $\mathcal{O}_N(1)$ has $Z$-weight $1$. \hfill $\Box$

\section{Construction}
Finally we present the calculations showing that under appropriate choices the notion of $\delta$-stability coincides with GIT-stability.
\begin{prop}\label{propconstr1}
Let $(a,q)$ be a point in $\bar{\mathcal{R}}$. For sufficiently large $l$ $(a,q)$ is (semi)stable in the GIT sense with respect to $\mathcal{O}_{\mathcal{Z}}(n_1,n_2)$ if and only if the following holds:
Set $W:=H^0(\mathcal{O}_X(l-m))$ and $q' := H^0(q(l))\colon V\otimes W \rightarrow H^0(\mathcal{E}(l))$. Then for every nontrivial proper subspace $U$ of $V$ we have: 
\begin{eqnarray}\label{eqconstr1}
\dim U[n_1P(l)-n_2]&(\leq)& P(m)[\dim(q'(U\otimes W))n_1 - \epsilon(U)n_2],
\end{eqnarray}
where $\epsilon(U) = 1$ if $U\subseteq \mathrm{im}a$ and $0$ otherwise.
\end{prop}
\em Proof: \em In order to apply the Hilbert$-$Mumford criterion we have to look at 1-parameter subgroups $\lambda\colon\mathbb{G}_m\rightarrow $SL$(V)$. Such a $\lambda$ is completely determined by giving a basis ${v_1,\dots,v_p}$ of $V$ and a weight vector $(\gamma_1,\dots,\gamma_p)\in\mathbb{Z}^p$ satisfying $\gamma_1 \leq \dots\leq\gamma_p$ and $\sum\gamma_i=0$. The action of $\lambda$ is then given by $\lambda(t)\cdot v_i = t^{\gamma_i}v_i$.\\
Now we look at a point $(a,q)\in  \bar{\mathcal{R}}$ represented by homomorphisms $q\colon V\otimes\mathcal{O}_X(-m) \rightarrow \mathcal{E}$ and $ a\colon H^0(\mathcal{D}(m))\rightarrow V$ and we denote by $\varphi\colon \mathcal{D}\rightarrow\mathcal{E}$ the corresponding framing. For the moment we fix an $l\geq m$ and let $W := H^0(\mathcal{O}_X(l-m))$ and $\varrho:=h^0(\mathcal{E}(l))=P(l)$. Now $q$ induces the homomorphisms $q'= H^0(q(l))$ and $q''\colon\Lambda^\varrho(V\otimes W)\rightarrow $ det$  H^0(\mathcal{E}(l))$. Let ${w_1,\dots,w_t}$ be a basis of $W$. We then get a basis of $\Lambda^\varrho(V\otimes W)$ by elements of the form
\[
u_{IJ}=(v_{i_1}\otimes w_{j_1})\wedge\dots\wedge(v_{i_\varrho}\otimes w_{j_\varrho})
\]
with multiindices $I$ and $J$ satisfying $i_\alpha\leq i_{\alpha+1}$ and $j_\alpha < j_{\alpha+1}$ if $i_\alpha= i_{\alpha+1}$. Now the action of $\lambda$ on $\Lambda^\varrho(V\otimes W)$ is given by
\begin{eqnarray*}
\lambda(t)\cdot u_{IJ} = t^{\gamma_I}u_{IJ}& with& \gamma_I:=\sum_\alpha\gamma_{i_\alpha}.
\end{eqnarray*}
Now $\mu(q'',\lambda)$ is given by  $-$min$\{\gamma_I|\exists I,J\text{ with }q''(u_{IJ})\not=0\}$. But a slightly better formulation is possible. We set $\psi(i)= \dim(q'(\langle v_1,\dots,v_i\rangle\otimes W))$. Then we have
\[\mu(q'',\lambda)=-\sum_{i=1}^p\gamma_i\big(\psi(i)-\psi(i-1)\big).\]
Observe that $\psi(i)-\psi(i-1)$ is always equal to one or zero. Thus the right hand side sums up all the $\gamma_i$ such that $\dim q'(\langle v_1,\dots,v_i\rangle\otimes W)$ is increasing. By the surjectivity of $q$ (or $q'$) we know that there are $\varrho$ such $\gamma_i$. Because the $\gamma_i$ are in increasing order this $\gamma_I$ must be the smallest such that $q'(\langle v_1,\dots,v_i\rangle\otimes W)=H^0(\mathcal{E}(l))$.\\
Next we want to determine $\mu(a,\lambda)$. Look at the following identification:
\begin{eqnarray*}
\text{Hom}(U,V)&\cong&U^\vee\otimes V\\
a: w_j\mapsto\sum_i\alpha_{ij}v_i&\leftrightarrow&\sum_{i,j}\alpha_{ij}w_j^\vee\otimes v_i.
\end{eqnarray*}
We now deduce easily that $\mu(a,\lambda)=\max\{\gamma_i|\exists j\text{ with }\alpha_{ij}\neq 0\}=\min\{\gamma_i|\mathrm{im}a\subseteq \langle v_1,\dots,v_i\rangle\}=\gamma_\tau$ where 
$\tau=\min\{i|\mathrm{im}a\subseteq \langle v_1,\dots,v_i\rangle\}$.\\
Now by the Hilbert$-$Mumford criterion $(q,a)$ is (semi)stable if and only if
\begin{eqnarray}\label{eqssroh}
n_1\cdot\mu(q'',\lambda)+n_2\cdot\mu(a,\lambda) (\geq)0,\text{ i.e.,}\\
n_1\cdot\sum_{i=1}^p\gamma_i\big(\psi(i)-\psi(i-1)\big) - n_2\cdot \gamma_\tau(\leq)0.\nonumber
\end{eqnarray}
Fixing a basis ${v_1,\dots,v_p}$ for the moment we can consider the left hand side as a linear form on the set of weight vectors. Thus it is enough to check the inequality for the special weight vectors
\begin{eqnarray*}
\gamma^{(i)}=(\underbrace{i-p,\dots,i-p}_i,\underbrace{i,\dots,i}_{p-i}),&i=1,\dots,p-1.
\end{eqnarray*}
Indeed, every weight vector can be expressed as a finite nonnegative linear combination of the $\gamma^{(i)}$. Now for such a $\gamma^{(i)}$ we have 
\begin{eqnarray*}
\gamma^{(i)}_\tau =
\begin{cases}i-p & \text{if im}a\subseteq\langle v_1,\dots,v_i\rangle\\
i&\text{otherwise}. \end{cases}
\end{eqnarray*}
In other words we have $\gamma^{(i)}_\tau = i-\epsilon(i)p$ where $\epsilon(i):=1$ if im$a\subseteq \langle v_1,\dots,v_i\rangle$ and $0$ otherwise.
On the other hand we have
\begin{eqnarray*}
\sum_{i=1}^p\gamma_i\big(\psi(i)-\psi(i-1)\big)  =  \varrho i-p\sum_{k=1}^i\big(\psi(k)-\psi(k-1)\big)  =  \varrho i-p\psi(i).
\end{eqnarray*}
Alltogether our inequality now reads:
\begin{eqnarray*}
i\cdot(n_1\varrho-n_2)&(\leq)&p\cdot(n_1\psi(i)-\epsilon(i)n_2).
\end{eqnarray*}
In particular this inequality does not contain no weights anymore. It does not even depend on the fixed basis we chose but on the subspaces spanned by this basis. Thus a point $(q,a)$ is (semi)stable if and only if for every nontrivial subspace $U\subseteq V$ we have:
\begin{eqnarray*}
\dim U\cdot (n_1\varrho-n_2)&(\leq)&\dim V\cdot\big(n_1\dim(q'(U\otimes W))-\epsilon(U)n_2\big)
\end{eqnarray*}
where $\epsilon(U):=1$ if im$a\subseteq U$ and $0$ otherwise.
But we have the identifications $\varrho=P(l)$ and $\dim V=P(m).$ \hfill $\Box$\\
\\
For every nontrivial subspace $U\subseteq V$ we denote by $\mathcal{F}_U$ the subsheaf of $\mathcal{E}$ generated by $U$.

\begin{lem}\label{lemgitinj}
For every GIT-semistable point $(a,q)$ the induced morphism $H^0(q(m))\colon V \rightarrow H^0(\mathcal{E}(m))$ is injective. In particular $\dim(V\cap H^0(\mathcal{E}(m)))\leq h^0(\mathcal{E}(m))$ where $V\cap H^0(\mathcal{E}(m))$ denotes the preimage of $H^0(\mathcal{E}(m))$ in $V$. Furthermore $q'$ is injective and for every subspace $U\subseteq V$ we have $\dim(q'(U\otimes W))\leq h^0(\mathcal{F}_U(l))$, where $\mathcal{F}_U:=q(U\otimes \mathcal{O}_X(-m))$. 
\end{lem}
\em Proof: \em Let $U\subseteq V$ denote the kernel of $H^0(q(m))$, then the generated subsheaf $\mathcal{F}_U$ is zero (and so is $\epsilon(\mathcal{F}_U)$), thus by the upper inequality $\dim U\leq 0$. Note that by the choice of $n_1,n_2$ we have $n_1\varrho-n_2\geq 0$.  \hfill $\Box$

\begin{prop}
For sufficiently large $l$ a point $(a,q)$ is GIT-(semi)stable if and only if for every nontrivial subspace $U\subseteq V$ we have the following inequality of polynomials in $l$:
\begin{eqnarray}\label{eqconstr2}
\dim U\cdot(n_1P(l)-n_2)&(\leq)& P(m)\cdot(n_1P_{\mathcal{F}_U}(l)-\epsilon(\mathcal{F}_U)n_2).
\end{eqnarray}
\end{prop}
\em Proof: \em
By Prop \ref{propconstr1} it is enough to show that (\ref{eqconstr1}) is equivalent to (\ref{eqconstr2}) for every $U$.\\
First of all we note that the family of these subsheaves $\mathcal{F}_U$ generated by some subspace of $U\subseteq V$ is bounded. Thus by taking $l$ large enough all the $\mathcal{F}_U$ are globally generated and we have $P_{\mathcal{F}_U}(l)=h^0(\mathcal{F}_U(l))=\dim(q'(U\otimes W))$ $\forall \mathcal{F}_U$.
We'll now have a look at the following diagram:
\[
\xymatrix{
H^0(\mathcal{D}(m))\otimes\mathcal{O}_X(-m) \ar[r]_/1em/a \ar[d]_{\mathrm{ev}} & V\otimes \mathcal{O}_X(-m) \ar[d]^q\\
\mathcal{D} \ar[r]^\varphi & \mathcal{E}.
}\]
It is easy to see, that $\mathrm{im}a\subseteq U \Rightarrow \mathrm{im}\varphi \subseteq \mathcal{F}_U$. Thus if $\epsilon(U)=1$, we have $\epsilon(\mathcal{F}_U)=1$ and in this case we easily see that (\ref{eqconstr1}) is equivalent to (\ref{eqconstr2}).\\
So the only interesting case is $\epsilon(U)=0$ and $\epsilon(\mathcal{F}_U)=1$.
Now let $(a,q)$ be GIT-(semi)stable with $\epsilon(U)=0$ and $\epsilon(\mathcal{F}_U)=1$. We claim that for every nontrivial subspace we also have inequality (\ref{eqconstr2}). Let $U$ be a nontrivial subspace of $V$ and let $U':=U\oplus \mathrm{im}a$. Note that $\epsilon(U')=1$ and $U'$ generates the same subsheaf $\mathcal{F}_U$. Thus we get the desired inequality with $\dim U$ replaced by $\dim U'$. But since $\dim U\leq\dim U'$ our claim follows.
Conversely let $(a,q)$ be a point satisfying (\ref{eqconstr2}) for every nontrivial subspace $U\subseteq V$. Now (\ref{eqconstr1}) is equivalent to $P(m)\cdot n_2 > 0$. But this is clear by the definition.\hfill $\Box$
 
\begin{prop}\label{propchl}
For sufficiently large $l$ a point $(a,q)$ is GIT-(semi)stable if and only if for every nontrivial subspace $U\subseteq V$ we have the following inequality of polynomials:
\begin{eqnarray}\label{eqconstr3}
 P\cdot\big(\dim U+\epsilon(\mathcal{F}_U)\delta(m)\big) +
 \delta\cdot\big(\dim U-\epsilon(\mathcal{F}_U)P(m)\big)& (\leq) &P_{\mathcal{F}_U}\cdot\big(P(m)+\delta(m)\big).\nonumber\\
\end{eqnarray}
\end{prop}
\em Proof: \em 
Again because the family of such subsheaves is bounded we can find an $l$ such that inequality (\ref{eqconstr2}) holds if and only if it holds as an inequality of polynomials in $l$. Now substitute
\[\frac{n_2}{n_1} = \frac{P(l)\cdot\delta(m)-\delta(l)\cdot P(m)}{P(m)+\delta(m)}.\]
We then easily derive the required inequality. \hfill $\Box$

\begin{thm}\label{thmcalc}
For sufficiently large $l$ if a point $(a,q)$ is GIT-semistable then the corresponding pair $(\mathcal{E},\varphi)$ is semistable (with respect to $\delta$) and $H^0(q(m))$ is an isomorphism. In particular every GIT-semistable point corresponds to a pair with pure sheaf.
\end{thm}
\em Proof: \em First we will drop the restriction to subsheaves generated by subspaces $U\subseteq V$. Let $\mathcal{F}$ be an arbitrary subsheaf of $\mathcal{E}$. Set $U:=V\cap H^0(\mathcal{F}(m))$. Then the subsheaf generated by $U$ is contained in $\mathcal{F}$. Now if $\epsilon(\mathcal{F}_U)=1$ we also have $\epsilon(\mathcal{F})=1$ and we easily get the required inequality. Conversely if $\epsilon(\mathcal{F})=1$ by definition of $U$ we have $\epsilon(U)=1$, hence $\epsilon(\mathcal{F}_U)=1$. Thus from now on we can consider arbitrary subsheaves together with the subspace $U:=V\cap H^0(\mathcal{F}(m))$.\\
Passing to the leading coefficient of the polynomials in the inequality (\ref{eqconstr3}) we have:
\begin{eqnarray*}
\dim U + \epsilon(\mathcal{F}) \delta(m) \leq \frac{r'}{r} (P(m)+\delta(m)).
\end{eqnarray*}
We now assign to any submodule $\mathcal{F}$ its corresponding quotient $\mathcal{G}$. Note that since $\dim U\leq h^0(\mathcal{F}(m))$ we find $\dim(V/U)\leq h^0(\mathcal{G}(m))$. We have
\begin{eqnarray*}
h^0(\mathcal{G}(m)) + \epsilon(\mathcal{G})\delta(m) & \geq & \dim(V/U) + \epsilon(\mathcal{G})\delta(m)
\\ & =  & \dim V + \delta(m) - (\dim U + \epsilon(\mathcal{F})\delta(m)) \\
& \geq & P(m) + \delta(m) - \frac{r'}{r} (P(m)+\delta(m)) \\
& \geq & \frac{r''}{r}(P(m)+\delta(m)).
\end{eqnarray*}
Since $(a,q)\in  \bar{\mathcal{R}}$, it deforms into a pair with torsion free sheaf. By Proposition \ref{propdeform} there is a torsion free sheaf $\mathcal{H}$ together with a homomorphism $\psi\colon\mathcal{E}\rightarrow\mathcal{H}$ satisfying $P(\mathcal{H})=P(\mathcal{E})$ and ker$(\psi)=\mathcal{T}(\mathcal{E})$. Next we want to show that ($\mathcal{H},\varphi_\mathcal{H})$ is semistable. As noted after Proposition \ref{propdeform} the corresponding homomorphism $\varphi_\mathcal{H}$ might be trivial. But here we will use the following lemma.

\begin{lem}
Let $l$ be chosen as in \em Proposition \ref{propchl} \em and for every GIT-semistable point $(a,q)$ let $\mathcal{T}$ denote the maximal subsheaf of strictly smaller dimension of $\mathcal{E}$ (cf.\ \cite{HL}, Definition 1.1.4). Then we have $\mathrm{im}\varphi\not\subseteq\mathcal{T}$, i.e., $\epsilon(\mathcal{T})=0$.

\end{lem}
\em Proof\em : Suppose $\epsilon(\mathcal{T})$ to be $1$ and look at the leading coefficients of (\ref{eqconstr3}). Since $\mathcal{T}$ is torsion the leading coefficient of $P_\mathcal{T}$ is zero. Hence we have
\begin{eqnarray*}
r\cdot (\dim U + \delta(m))&\leq&0
\end{eqnarray*}
But since $\delta(m)$ is clearly positive we have a contradiction. \hfill $\Box$\\
\\
\em Continuation of the proof of \em Theorem \ref{thmcalc}: Thus we have seen that $\varphi_\mathcal{H}$ is nontrivial and we can now continue the proof by showing that $\mathcal{H}$ is semistable. If $\pi\colon\mathcal{H}\rightarrow \mathcal{G}_\mathcal{H}$ is any quotient of $\mathcal{H}$, let $\mathcal{G}$ denote the image of $\mathcal{E}$ by $\pi\circ\psi$. This is a quotient of $\mathcal{E}$ which is contained in $\mathcal{G}_\mathcal{H}$. Hence we have $h^0(\mathcal{G}_\mathcal{H}(m))\geq h^0(\mathcal{G}(m))$ and $\epsilon(\mathcal{G}_\mathcal{H})=0 \Rightarrow \epsilon(\mathcal{G})=0.$ And even if $\epsilon(\mathcal{G}_\mathcal{H})=1$ and $\epsilon(\mathcal{G})=0$ we get the following inequalities:
\begin{eqnarray}\label{eqconstr4}
h^0(\mathcal{G}_\mathcal{H}(m)) + \epsilon(\mathcal{G}_\mathcal{H})\delta(m)
& \geq & h^0(\mathcal{G}(m)) + \epsilon(\mathcal{G})\delta(m) \\
&\geq& \frac{r''}{r}(P(m)+\delta(m)) = \frac{r_{\mathcal{G}_\mathcal{H}}}{r}(P(m)+\delta(m))\nonumber.
\end{eqnarray}
By Proposition \ref{propssrephrase} $(\mathcal{H},\varphi_\mathcal{H})$ is semistable and therefore $m$-regular. By taking $\mathcal{G}_\mathcal{H}=\mathcal{H}$ we find in fact equality in (\ref{eqconstr4}) everywhere. This shows, that $h^0(\psi(\mathcal{E})(m))=h^0(\mathcal{H}(m))=P(m)$ and since $\mathcal{H}$ is globally generated we find $\psi$ to be surjective. Since $\mathcal{E}$ and $\mathcal{H}$ have the same Hilbert polynomial $\psi$ is in fact an isomorphism. Hence $\mathcal{E}$ is (semi)stable.\\
Because of Lemma \ref{lemgitinj} $H^0(q(m))$ is injective and surjectivity now follows easily from dimension reasons. \hfill $\Box$

\begin{thm}
Let $(\mathcal{E},\varphi)$ be a $\delta$-(semi)stable pair such that $q$ induces an isomorphism $V\rightarrow H^0(\mathcal{E}(m))$. Then the corresponding point $(a,q)$ is (semi)stable in the GIT-sense.
\end{thm}
\em Proof: \em By Proposition \ref{propssrephrase} we have
\begin{eqnarray*}
h^0(\mathcal{F}(m)) + \epsilon(\mathcal{F})\delta(m) &(\leq) &\frac{r'}{r}(P(m)+\delta(m))
\end{eqnarray*}
for every subsheaf $\mathcal{F}$ of $\mathcal{E}$ of rank $r'$ satisfying $0<r'<r=rk(\mathcal{E})$. 
If $(\mathcal{E},\varphi)$ is stable then the inequality is strict. Now if $U$ is an arbitrary subspace of $V$ and $\mathcal{F}_U$ denotes the subsheaf generated by $U$ we have $U\subseteq V\cap H^0(\mathcal{F}_U(m))$ and $\dim U \leq h^0(\mathcal{F}_U(m))$. Thus we get the following strict inequality:
\begin{eqnarray*}
\dim U + \epsilon(\mathcal{F}_U)\delta(m) &< &\frac{r'}{r}(P(m)+\delta(m)).
\end{eqnarray*}
This is a strict inequality of the leading coefficients of the desired inequality.\\
If $(\mathcal{E},\varphi)$ is semistable but not stable, again strict inequality holds except the case of a destabilizing semistable subsheaf $\mathcal{F}$. But such an $\mathcal{F}$ is semistable with the same reduced Hilbert polynomial. Hence $\mathcal{F}$ has the same slope as $\mathcal{E}$ and by the choice of $m$ we have $P_\mathcal{F}(m)=h^0(\mathcal{F}(m))$. The destablizing condition on $\mathcal{F}$ says:
\[P_\mathcal{F} + \epsilon(\mathcal{F})\delta=\frac{r'}{r}(P+\delta).\]
Now let $U:=V\cap H^0(\mathcal{F}(m))$ and note that $\dim U=h^0(\mathcal{F}(m))=P_\mathcal{F}(m)$ because $q$ induces an isomorphism $V\rightarrow H^0(\mathcal{E}(m))$. The terms (let's call them (1) and (2)) we want to show to be equal are:
\begin{eqnarray*}
(1) & P_\mathcal{F}(P(m)+\delta(m)) =\big(\frac{r'}{r}(P+\delta\big)-\epsilon(\mathcal{F})\delta\big)\big(P(m)+\delta(m)\big)=\\ 
& P\cdot\frac{r'}{r}\big(P(m)+\delta(m)\big)+
\delta\cdot\big(\frac{r'}{r}(P(m)+\delta(m))-\epsilon(\mathcal{F})(P(m)+\delta(m))\big),\\\\
(2) & P\cdot\big(\dim U+\epsilon(\mathcal{F})\delta(m)\big) + \delta\cdot\big(\dim U-\epsilon(\mathcal{F})P(m)\big)=\\
& P\cdot\big(P_\mathcal{F}(m)+\epsilon(\mathcal{F})\delta(m)\big) +\delta\cdot\big(P_\mathcal{F}(m)-\epsilon(\mathcal{F})P(m)\big).
\end{eqnarray*}
By applying the destabilizing condition evaluated at $m$ to the coefficients standing infront of $P$ and $\delta$ shows that they are indeed the same. Thus the point $(a,q)$ corresponding to $(\mathcal{E},\varphi)$ is semistable but not stable. \hfill $\Box$

\section{Variation of the Stability Parameter}
In this section we want to study how the moduli space $M_{X,\delta}(\mathcal{D},P)$ changes if we vary the stability parameter $\delta$. Therefore we fix $P$ and $\mathcal{D}$.

\begin{lem}
Define $\mathcal{W}$ to be the set of all subsheaves \em im\em$\varphi$ occuring as the image of a homomorphism $\varphi$ in a $\delta$-semistable pair for any $\delta$. Then $\mathcal{W}$ is bounded.
\end{lem}
\em Proof: \em For every semistable pair $(\mathcal{E},\varphi)$ there is a corresponding point $(a,q)\in N\times Q$ (cf.\ Lemma \ref{lemassaq}) such that 
$q\circ a = \varphi\circ \mathrm{ev}$. On $N\times Q\times X$ there is the universal homomorphism (cf.\ Proposition \ref{propclsub})
\[\tilde{q}\circ\tilde{a}\colon H^0(\mathcal{D}(m))\otimes\pi^\star_X\mathcal{O}_X(-m)\rightarrow \pi^\star_Q \tilde{\mathcal{E}}. \]
It has the universal property that its fibre at the point $(a,q)\in N\times Q$ (which is a homomorphism of sheaves on $X$) is just $q\circ a$. Thus we find that every sheaf in $\mathcal{W}$ occurs as a fibre of the sheaf im$(\tilde{q}\circ\tilde{a})$. Hence $\mathcal{W}$ is bounded.\hfill $\Box$\\
\\
The next result is a refinement of Proposition \ref{propgensur}.

\begin{prop}
There is a rational polynomial $\delta_{\max}$ of degree \textbf{\textit{less}} than $\dim X$ such that for every $\delta>\delta_{\max}$ and every pair $(\mathcal{E},\varphi)$ the following two assumptions are equivalent:
\begin{itemize}
\item[\em(i)\em] $(\mathcal{E},\varphi)$ is $\delta$-semistable,
\item[\em(ii)\em] $\varphi$ is generically surjective.
\end{itemize}
\end{prop}
\em Proof: \em We follow closely the proof of Proposition \ref{propgensur}:\\
(i) $\Rightarrow$ (ii): If $\varphi$ was not generically surjective then the saturation $\mathcal{F}$ of im$\varphi$ would be a proper subsheaf of $\mathcal{E}$. Semistability yields
\begin{eqnarray}
\label{eqgensur} p_\mathcal{F}\leq p_\mathcal{E} + \delta\cdot\underbrace{(\frac{1}{r_\mathcal{E}}-\frac{1}{r_\mathcal{F}})}_{=:c}.
\end{eqnarray}
Now one can easily see that $-1<c<0$, because $1\leq r_{\mathcal{F}}\leq r_{\mathcal{E}}$. Since the leading coefficients of $p_\mathcal{F}$ and $p_\mathcal{E}$ agree we may choose a $\delta$ of degree less than $\dim X$ such that the inequality (\ref{eqgensur}) is violated. Since the set of image sheaves of semistable pairs is bounded there are only finitely many polynomials we have to consider. Hence we may choose a $\delta_{\max}$ working for all of these.\\
(ii) $\Rightarrow$ (i): Conversely, let $(\mathcal{E},\varphi)$ be a pair, where $\varphi$ is generically surjective. Then for any subsheaf $\mathcal{F}\subseteq \mathcal{E}$ we have $\epsilon(\mathcal{F})=0$. Thus if $(\mathcal{E},\varphi)$ is not $\delta_{\max}$-semistable then $\mathcal{F}$ satisfies
\begin{eqnarray}\label{eqchambstr}
\mu_{\mathcal{F}}  >  \mu_{\mathcal{E}} + \frac{\delta_{\max,1}}{r_\mathcal{E}} >  \mu_{\mathcal{E}},
\end{eqnarray} 
where $\delta_{\max,1}$ shall denote the coefficient of $\delta_{\max}$ in degree $d-1$. Note that the set of sheaves $\mathcal{E}$ occuring in a pair $(\mathcal{E},\varphi)$ where $\varphi$ is generically surjective is certainly bounded. Therefore the set of subsheaves satisfying (\ref{eqchambstr}) is bounded and we may choose $\delta_{\max,1}$ big enough to contradict the first inequality in (\ref{eqchambstr}) for all $\mathcal{F}$. \hfill $\Box$

\begin{lem}
For any $\delta\geq 0$ we consider the set $\mathcal{S}_{\delta}$ of all subsheaves occuring as a $\delta$-destabilizing subsheaf of any pair $(\mathcal{E},\varphi)$ which is $\delta'$-semistable for some $\delta'$. Let $\mathcal{S}:= \bigcup_{\delta\geq 0} \mathcal{S}_{\delta}$. Then $\mathcal{S}$ is bounded.
\end{lem}
\em Proof: \em If $\delta\leq\delta_{\max}$ then the destabilizing condition reads
\[ \mu_\mathcal{F} > \mu_\mathcal{E} + \delta_1(\frac{1}{r_\mathcal{E}}-\frac{\epsilon(\mathcal{F})}{r_\mathcal{F}}).\]
But the right hand side is bounded from below by $\mu_\mathcal{E}-\delta_{\max,1}$. Thus for every sheaf in $\mathcal{S}$ there is a uniform bound $\mu_\mathcal{F}> \mu_\mathcal{E}-\delta_{\max,1}$. If $\delta>\delta_{\max}$ then by the proposition above there are no destabilizing subsheaves. Now the claim follows from \cite{Gro}, Lemme 2.5.\hfill $\Box$

\begin{cor}
There are only finitely many polynomials occuring in a destabilizing condition for a semistable pair in $M_{X,\delta}(\mathcal{D},P)$ for any $\delta$.
\end{cor}
From these results one may deduce easily the following theorem which summarizes how stability depends on the parameter $\delta$. 

\begin{thm}
There are finitely many critical values $\delta^1,\dots,\delta^s\in \mathbb{Q}[z]$,
\[\xymatrix@R=0pt@C=1pt{  \ar  '[rrrr] [rrrrrrrrrrrr]  & | & & | & \cdots & & |  &&&&&&\\
& 0 && \delta^1 &&& \delta^s &&&&&&\infty
}\]
such that setting $\delta^0:=0$ and $\delta^{s+1}:=\infty$ the following properties hold true:
\begin{itemize}
\item[\em(i)\em] For $i=0,\dots,s$ and $\delta,\delta '\in(\delta^i,\delta^{i+1})$, one has
\[\mathcal{R}_\delta^{(s)s}=\mathcal{R}_{\delta '}^{(s)s}.\]
\item[\em(ii)\em] For $i=0,\dots,s$ and $\delta\in(\delta^i,\delta^{i+1})$, there are the inclusions
\begin{eqnarray*}
\mathcal{R}_\delta^{ss}\subseteq \mathcal{R}_{\delta^i}^{ss}\cap \mathcal{R}_{\delta^{i+1}}^{ss},\\
\mathcal{R}_\delta^s\supseteq \mathcal{R}_{\delta^i}^s\cup \mathcal{R}_{\delta^{i+1}}^s.
\end{eqnarray*}
\item [\em(iii)\em] For $i=0,\dots,s$ and $\delta\in(\delta^i,\delta^{i+1})$, one has
\[ M_{X,\delta}(\mathcal{D},P)=M^s_{X,\delta}(\mathcal{D},P).\]
\end{itemize}

\end{thm}
\em Remark: \em Of course we have $\delta^s=\delta_{\max}$. Alltogether these inclusions yield a nice diagram often called the \em chamber structure \em of the stability parameter:
\[\xymatrix@C=1pt{ & M^0_X(\mathcal{D},P) \ar@{.}[rrrr] \ar[dl] \ar[dr] && & & & M^s_X(\mathcal{D},P)\ar[dl] \ar[dr]\\
M_{X,\delta^0}(\mathcal{D},P) & & M_{X,\delta^1}(\mathcal{D},P)\ar@{.}[rrr] &&&  M_{X,\delta^s}(\mathcal{D},P) & &  M_{X,\delta^{s+1}}(\mathcal{D},P),}\]
where we set $M^i_X(\mathcal{D},P):=M_{X,\delta}(\mathcal{D},P)$, for some $\delta\in(\delta^i,\delta^{i+1}), i=0,\dots,s$. The diagonal maps are easily obtained by the universal property of the categorical quotient.

\section{Coherent Systems}
\begin{defi}
A \em coherent system \em is a pair $(\Gamma,\mathcal{E})$ consisting of a coherent sheaf $\mathcal{E}$ on $X$ and a vector subspace $\Gamma\subseteq H^0(\mathcal{E})$ of dimension $r$. Two coherent systems $(\Gamma,\mathcal{E})$ and $(\Gamma ',\mathcal{E}')$ are called \em isomorphic \em if there is an isomorphism of sheaves $\mathcal{E}\rightarrow \mathcal{E}'$ such that the corresponding map of global sections maps $\Gamma$ isomorphically onto $\Gamma'$.
\end{defi}
\em Remark: \em As a simple consequence we deduce that for isomorphic systems we have $\dim\Gamma=\dim\Gamma'$.\\
\\
Next we want to introduce a notion of stability for coherent systems. Just like the stability of pairs it depends on a parameter, a nonnegative rational polynomial $\alpha$.

\begin{defi}
Let $(\Gamma,\mathcal{E})$ be a coherent system. For a subsheaf $\mathcal{F}\subseteq \mathcal{E}$ we define $\Gamma ':=\Gamma\cap H^0(\mathcal{F})$.
A coherent system $(\Gamma,\mathcal{E})$ is called \em $\alpha$-(semi)stable \em if for every nontrivial saturated subsheaf $\mathcal{F}\subsetneq\mathcal{E}$ we have the following inequality of rational polynomials:
\begin{eqnarray*}
\dim\Gamma '\cdot \alpha + P_\mathcal{F}
&(\leq)& \frac{r_\mathcal{F}}{r_\mathcal{E}}\big{\{}\dim\Gamma\cdot \alpha + P_\mathcal{E}\big{\}}.
\end{eqnarray*}
\end{defi}
\em Remark: \em Consider the case $\dim\Gamma=0$. Then our stability condition reduces to the usual one for sheaves. Thus from now on we may assume $\dim\Gamma>0$.\\
\\
The natural evaluation map 
\[ ev_{(\Gamma,\mathcal{E})}\colon \Gamma\otimes\mathcal{O}_X\rightarrow \mathcal{E}\]
provides a pair $(\mathcal{E},ev_{(\Gamma,\mathcal{E})})$ which we want to call \em the corresponding pair to the system $(\Gamma,\mathcal{E})$\em . Note that the evaluation map, of course, is injective on global sections.\\
Conversely any pair $(\mathcal{E},\varphi)$ consisting of a coherent sheaf $\mathcal{E}$ and a map $\varphi\colon\mathcal{O}_X^r\rightarrow \mathcal{E}$ yields a coherent system $($im$(H^0(\varphi)),\mathcal{E})$. If in addition $H^0(\varphi)$ is injective the dimension of im$(H^0(\varphi))$ is, of course, equal to $r$.\\
\\
Just as we have done in Sections 2-4 in the case of $\delta$-(semi)stable  pairs one can define a moduli functor for $\alpha$-(semi)stable coherent systems and construct the moduli space $\mathrm{Syst}_{X,\alpha}(r,P)$ of $\alpha$-semistable coherent systems $(\Gamma,\mathcal{E})$, where $P$ is the Hilbert polynomial of $\mathcal{E}$ and $\dim\Gamma=r$. This construction is done by Le Potier in \cite{LeP} and \cite{LeP2}. In the following passage we want to show that $\mathrm{Syst}_{X,\alpha}(r,P)$ can be obtained as a quotient by some GL$_r$-action of a moduli space of semistable pairs which we constructed in Section 4. We need the following results:

\begin{lem}\label{lemsystcorr}
Let $(\Gamma,\mathcal{E})$ be a $\alpha$-(semi)stable coherent system. Set $\delta:=\dim\Gamma \cdot \alpha$. Then the corresponding pair $(\mathcal{E},ev_{(\Gamma,\mathcal{E})})$ is $\delta$-(semi)stable.
\end{lem}
\em Proof: \em Let $\mathcal{F}\subsetneq\mathcal{E}$ be a saturated subsheaf. Assume $\mathrm{im}(\mathrm{ev}_{(\Gamma,\mathcal{E})})\subseteq\mathcal{F}$. Thus $\Gamma\subseteq H^0(\mathcal{F})$ and we have $\Gamma '=\Gamma\cap H^0(\mathcal{F})=\Gamma$. Now $(\Gamma,\mathcal{E})$ is $\alpha$-(semi)stable and we can apply the stability condition to $\mathcal{F}$:
\begin{eqnarray*}
\dim\Gamma \cdot \alpha + P_\mathcal{F} &(\leq) &\frac{r_\mathcal{F}}{r_\mathcal{E}}\big{\{}\dim\Gamma\cdot \alpha + P_\mathcal{E}\big{\}}.
\end{eqnarray*}
Substituting $\delta=\dim\Gamma \cdot \alpha$ we get the required inequality for $\mathcal{F}$.\\
Now assume $\mathrm{im}(\mathrm{ev}_{(\Gamma,\mathcal{E})})\nsubseteq\mathcal{F}$. Since $\dim\Gamma '$ and $\alpha$ are nonnegative the required inequality this time follows immediately:
\begin{eqnarray*}
P_\mathcal{F} \leq \dim\Gamma '\cdot \alpha + P_\mathcal{F} &(\leq)& \frac{r_\mathcal{F}}{r_\mathcal{E}}\big{\{}\dim\Gamma\cdot \alpha + P_\mathcal{E}\big{\}}.
\end{eqnarray*}
\hfill $\Box$\\
\\
Set $\mathcal{D}=\mathcal{O}_X^r$ and consider the set $\mathcal{R}\subseteq Q\times N$ of Definition/Lemma \ref{defilemr}. Now Lemma \ref{lemsystcorr} above says that all semistable coherent systems of type $(r,P)$ are parametrized by $\mathcal{R}$. Consider the $\mathrm{SL}_r(k)\times \mathrm{SL}(V)$-action on $\mathcal{R}$. Again there is a natural linearization of this action in $\mathcal{O}_\mathcal{R}(n_1,n_2)$ where we choose $n_1,n_2$ as before. Analogously to Proposition \ref{propconstr1} we calculate the weights for this linearization. 
\begin{prop}\label{propchss}
For every GIT-semistable point $(a,q)\in\mathcal{R}$ the corresponding coherent system $(\Gamma,\mathcal{E})$ is semistable.
\end{prop}
\em Proof: \em Let $\mathcal{F}\subseteq\mathcal{E}$, $\Gamma ':=\Gamma\cap H^0(\mathcal{F})$ and $j:=\dim\Gamma '$. Fix a basis $x_1,\dots,x_j$ of $\Gamma '$ and extend it to a basis $x_1,\dots,x_r$ of $\Gamma$. Furthermore set $U:= V\cap H^0(\mathcal{F}(m))$ and $i:=\dim U$. Clearly we have $a(\Gamma '\otimes H^0(\mathcal{O}_X(m)))\subseteq U$. We divide into two cases.\\
Case \textcircled{1}: im$a\nsubseteq U$. Choose a basis $v_1,\dots,v_i$ of $U$ and extend it to a basis $v_1,\dots,v_p$ of $V$. Consider the weight vectors
\[ \frac{p}{r}\cdot(\underbrace{-j,\dots,-j}_{r-j},\underbrace{r-j,\dots,r-j}_j) \text{ and } (\underbrace{i-p,\dots,i-p}_i,\underbrace{i,\dots,i}_{p-i}).\]
Let $\lambda$ be the formal one-parameter subgroup associated to these vectors. One easily derives
\begin{eqnarray}\label{eqweight1}
 \mu(a,\lambda) = i - j\cdot \frac{p}{r}.
\end{eqnarray}
From the semistability condition on the point $(a,q)$ we derive as in (\ref{eqssroh})
\begin{eqnarray}\label{eqchss}
-i\cdot P(l) + p\cdot P_\mathcal{F}(l) + \mu(a)\frac{\delta(m)\cdot P(l)-p\cdot\delta(l)}{p+\delta(m)}&\geq& 0
\end{eqnarray}
and the leading coefficient on the left hand side is (up to a factor)
\begin{eqnarray*}\label{eqleadco}
\frac{\delta(m)}{p}(\mu(a)-i)-i+\frac{r_\mathcal{F}}{r_\mathcal{E}}(p+\delta(m)).
\end{eqnarray*}
We plug in equation (\ref{eqweight1}) and end up with the desired semistability condition on $\mathcal{F}$.\\
Case \textcircled{2}: im$a\subseteq U$. Now we have $\Gamma '=\Gamma$, thus $r=j$. Therefore we choose the weight vectors
\[ (0,\dots,0) \text{ and } (\underbrace{i-p,\dots,i-p}_i,\underbrace{i,\dots,i}_{p-i})\]
and find
\begin{eqnarray*}\label{eqweight}
 \mu(a,\lambda)=i-p.
\end{eqnarray*}
Proceeding just as in case \textcircled{1} we again derive the necessary semistability condition. \hfill $\Box$\\
\\
\begin{prop}
Consider a point $(a,q)\in\mathcal{R}$ such that the corresponding coherent system is semistable. Then $(a,q)$ is GIT-semistable.
\end{prop}
\em Proof: \em We now have to consider all the one parameter subgroups of $\mathrm{SL}_r(k)\times \mathrm{SL}(V)$ and show that the corresponding weight is nonnegative. Since it is nearly impossible to study all kinds of such subgroups it is more convenient to restrict to special ones. Thus we apply a result by A.\ Schmitt (cf.\ \cite{Schm2}, Theorem 3.4) allowing us to restrict to one parameter subgroups corresponding to weight vectors of the special form 
\begin{eqnarray*}
\big{(}\frac{1}{r}(j-r,\dots,j-r,j,\dots,j);\frac{1}{p}(i-p,\dots,i-p,i,\dots,i)\big{)},
\end{eqnarray*}
such that the corresponding weight is exactly $i-j\cdot p/r$. Now we are almost done. Just as in Proposition \ref{propchss} we see that the weight is nonnegative if and only if (\ref{eqchss}) is fulfilled. Therefore we first look at the leading coefficient of the left hand side of (\ref{eqchss}) wich is nonnegative if and only if
\begin{eqnarray}
j\cdot\alpha + P_\mathcal{F}&\leq & \frac{r_\mathcal{F}}{r_\mathcal{E}}(P_\mathcal{E}+r\cdot\alpha),
\end{eqnarray}
where $\mathcal{F}$ is the subsheaf of $\mathcal{E}$ generated by the subspace spanned by the first $i$ elements of $V$ corresponding to the weight vector. Now since $j\leq \dim \Gamma '$ this is a strict inequality unless $(\Gamma',\mathcal{F})$ is destabilizing. In this case we proceed just as we did at the ending of Section 4 observing that in this case $j=\dim\Gamma'$.$\hfill \Box$
\\
\\
As the last ingredient we now state a powerful tool for the construction of quotients if we have a group action of a product of groups.

\begin{thm}\label{thmprodact}
Let $X$ be a projective $k$-scheme with group actions $\sigma\colon G\times X \rightarrow X$ and $\tau\colon H\times X \rightarrow X$ of reductive groups $G$ and $H$ coming with linearizations $\bar{\sigma}$ and $\bar{\tau}$. Suppose these group actions do commute such that there is an action of the product $\sigma\times \tau\colon G\times H\times X \rightarrow X$. Set $Q_\sigma:= X\sslash_{\bar{\sigma}} G$ and let $\pi\colon X^{ss}_{\bar{\sigma}}\rightarrow Q_\sigma$ be the quotient map. Then there is a natural $H$-action $\gamma$ on $Q_\sigma$ together with a linearization $\bar{\gamma}$ such that the set of $\bar{\sigma}\times \bar{\tau}$-semistable points is $\pi^{-1}((Q_\sigma)_\gamma^{ss})$ and 
\[X\sslash_{\bar{\sigma}\times \bar{\tau}}(G\times H) \cong Q_\sigma\sslash_{\bar{\gamma}} H\]
\end{thm}
\em Proof: \em \cite{Schm}, Theorem 1.5.3.1. \hfill $\Box$\\
\\
Alltogether we have proven the following result:

\begin{thm}
There is a rational map $M_{X,r\cdot\alpha}(\mathcal{O}_X^r,P) \dashrightarrow Syst_{X,\alpha}(r,P)$ which is a quotient map of the natural $\mathrm{GL}_r$-action on $M_{X,r\cdot\alpha}(\mathcal{O}_X^r,P).$ 
\end{thm}

\end{document}